\documentclass[12pt]{article}
\usepackage{mathrsfs}
\usepackage{stmaryrd}
\usepackage{amsfonts,amsmath,amssymb,amscd}
\usepackage{shadow}
\usepackage{graphicx,floatrow}
\usepackage{color}
\usepackage{longtable}
\usepackage{pstricks,multido}
\usepackage{hyperref}
\usepackage{qtree}
\usepackage{tabularx}
\usepackage{array}
\usepackage{caption}
\allowdisplaybreaks

\newtheorem{thm}{Theorem}[section]
\newtheorem{lem}[thm]{Lemma}
\newtheorem{prop}[thm]{Proposition}

\newtheorem{conj}[thm]{Conjecture}

\def\qed{\hfill \rule{4pt}{7pt}}

\def\pf{\noindent {\it{Proof.} \hskip 2pt}}

\numberwithin{equation}{section}

\begin{document}
\begin{center}
{\large\bf   A new bijective proof of Babson and Steingr\'{\i}msson's conjecture}
\end{center}

\begin{center}
Joanna N. Chen$^1$, Shouxiao Li$^2$

$^{1}$College of Science\\
Tianjin University of Technology\\
Tianjin 300384, P.R. China

$^2$
College of Computer and Information Engineering\\
Tianjin Agricultural University\\
Tianjin 300384, P.R. China

$^1$joannachen@tjut.edu.cn,
 $^2$shouxiao09009@163.com.

\end{center}

\begin{abstract}
Babson and Steingr\'{\i}msson introduced generalized permutation patterns and showed that most of the Mahonian statistics in the literature can be expressed by
the combination of generalized pattern functions. Particularly, they defined
a new Mahonian statistic in terms of generalized pattern functions, which is denoted  $stat$.
Given a permutation $\pi$, let $des(\pi)$  denote  the
 descent number of $\pi$ and $maj(\pi)$ denote
 the major index of  $\pi$.
Babson and Steingr\'{\i}msson conjectured that
$(des,stat)$ and $(des,maj)$ are equidistributed on $S_n$.
Foata and Zeilberger settled this conjecture using
q-enumeration, generating functions and Maple packages ROTA and PERCY.
Later, Burstein provided a bijective proof of a refinement of this conjecture.
In this paper, we give a new bijective proof of this conjecture.
\end{abstract}

\noindent {\bf Keywords}: Euler-Mahonian, bijection, involution

\noindent {\bf AMS  Subject Classifications}: 05A05, 05A15

\section{Introduction}
In this paper, we give a new bijective proof of a conjecture of  Babson and Steingr\'{\i}msson \cite{Babson} on Euler-Mahonian statistics.

Let $S_n$ denote the set of all the permutations of $[n]=\{1,2,\cdots,n\}$.
Given a permutation $\pi=\pi_1 \pi_2 \cdots \pi_n \in S_n$,
a descent of $\pi$ is a position $i \in [n-1]$ such that $\pi_i > \pi_{i+1}$, where $\pi_i$ and $\pi_{i+1}$ are called a descent top and a descent bottom, respectively.
An ascent of $\pi$ is a position $i \in [n-1]$ such that
$\pi_i < \pi_{i+1}$, where $\pi_i$ is called an ascent
bottom and $\pi_{i+1}$ is called an ascent top.
The descent set and the ascent set of $\pi$
are given by
\[
Des(\pi)=\{i\colon \pi_i>\pi_{i+1}\},\]
\[Asc(\pi)=\{i\colon \pi_i<\pi_{i+1}\}.
\]
The set of the inversions of $\pi$ is
$$
Inv(\pi)=\{(i,j)\colon 1\leq i< j \leq n, \pi_i > \pi_j\}.
$$
Let $des(\pi)$, $asc(\pi)$ and $inv(\pi)$
be the descent number, the ascent number and the inversion number of $\pi$, which are defined by
$des(\pi)=|Des(\pi)|$,
 $asc(\pi)=|Asc(\pi)|$ and $inv(\pi)=|Inv(\pi)|$, respectively.
The major index of  $\pi$, denoted $maj(\pi)$, is given by
\[maj(\pi)=\sum_{i \in  Des(\pi)}  i.\]

Suppose that $st_1$ is a statistic on the object $Obj_1$ and $st_2$ is a statistic on the object $Obj_2$. If
\[
\sum_{\sigma \in Obj_1} q^{st_1(\sigma)}=\sum_{\sigma \in Obj_2} q^{st_2(\sigma)},
\]
we say that the statistic $st_1$ over $Obj_1$ is equidistributed with the statistic $st_2$ over $Obj_2$.

A statistic on $S_n$ is said to be Eulerian if it is
equidistributed with the statistic $des$ on $S_n$.
While a statistic on $S_n$ is said to be Mahonian if
it is equidistributed with the statistic  $inv$ on $S_n$. It is well-known that
\[\sum_{\pi \in S_n} q^{inv(\pi)}=\sum_{\pi \in S_n} q^{maj(\pi)}=[n]_q !,\]
where $[n]_q=1+q+\cdots +q^{n-1}$ and $[n]_q!=[n]_q [n-1]_q \cdots [1]_q$. Thus, the major index $maj$ is a Mahonian statistic. A pair of statistics on $S_n$ is said to be Euler-Mahonian if it is equidistributed with the joint distribution of the descent number and the major index.

In \cite{Babson}, Babson and Steingr\'{\i}msson introduced generalized permutation patterns, where they allow the
requirement that two adjacent letters in a pattern must be
adjacent in the permutation.
Let $\mathcal{A}$ be the alphabet $\{a,b,c,\ldots\}$ with the usual ordering. We write patterns as words in $\mathcal{A}$, where two adjacent letters
may or may not be separated by a dash. Two adjacent letters without a dash in a pattern indicates that the corresponding letters in the permutation must be
adjacent. Given a generalized pattern $\tau$ and a permutation $\pi$, we say a subsequence of $\pi$ is an occurrence (or instance) of $\tau$ in $\pi$ if it is order-isomorphic to $\tau$ and satisfies the above dash conditions.
Let $(\tau)\pi$ denote
 the number of occurrences of $\tau$ in $\pi$. Here, we see $(\tau)$ as a
  generalized pattern function.
For example, an occurrence of the generalized pattern $b$-$ca$ in a permutation
$\pi=\pi_1\pi_2 \cdots \pi_n$ is a subsequence $\pi_i\pi_j\pi_{j+1}$ such that
$i<j$ and $\pi_{j+1}<\pi_i<\pi_j$.
For $\pi=4753162$, we have $(b$-$ca)\pi=4$.

Further,  Babson and Steingr\'{\i}msson \cite{Babson} showed that almost all of
the Mahonian permutation statistics in the literature can be written as linear
combinations of generalized patterns. We list some of them below.
\[maj=(a-cb)+(b-ca)+(c-ba)+(ba),\]
\[stat=(ac-b)+(ba-c)+(cb-a)+(ba).\]

They  conjectured that the statistic $(des,stat)$ is Euler-Mahonian.
\begin{conj}\label{stat}
The distribution of the bistatistic  $(des,stat)$ is equal to that of $(des,maj)$.
\end{conj}

In 2001, D. Foata and  D. Zeilberger \cite{Foata} gave a proof of this conjecture  using
q-enumeration and generating functions and an almost completely automated proof via Maple packages ROTA and PERCY.

Given a permutation $\pi=\pi_1 \pi_2 \cdots \pi_n$, let $F(\pi)=\pi_1$ be the
first letter of $\pi$ and
\[adj(\pi)=|\{i \colon 1 \leq i \leq n \text{~and~} \pi_i- \pi_{i+1}=1\}|,\]
where $\pi_{n+1}=0$.
 Burstein \cite{Burstein} provided a
bijective proof of the following  refinement of Conjecture \ref{stat} as follows.
\begin{thm}
Statistics $(adj, des, F, maj, stat)$ and $(adj, des, F, stat, maj)$
are equidistributed over $S_n$ for all $n$.
\end{thm}

In this paper, we will give a new bijective proof of Conjecture \ref{stat}, which does not preserve the statistic $adj$.

\section{A new bijective proof of Conjecture \ref{stat}}

In this section, we recall a particular bijection $\varphi$ on $S_n$ that
maps the inversion number to the major index, which is due to Carlitz \cite{Carlitz} and stated more clearly in \cite{remmel} and \cite{skandera}. Based on this, we give an analogous bijection which proves Conjecture \ref{stat}.

To give a description of $\varphi$, we first recall two labeling schemes for  permutations. It involves accounting for the effects of inserting a new largest element into a permutation, and so it is known as the insertion method.

Given a permutation $\sigma=\sigma_1 \sigma_2 \cdots \sigma_{n-1}$
in $S_{n-1}$. To obtain a permutation
$\pi \in S_n$, we can insert $n$ in $n$ spaces, namely,
immediately before $\sigma_1$ or  immediately after $\sigma_i$
for $1 \leq i \leq n$. In order to keep track of how the
insertion of $n$ affects the inversion number and the major index,
we may define two labelings of the $n$ inserting spaces.

The inv-labeling of $\sigma$ is given by numbering the spaces from
right to left with $0,1,\ldots,n-1$. The maj-labeling of $\sigma$ is
obtained by labeling the space after $\sigma_{n-1}$ with $0$, labeling
the descents from right to left with $1,2,\ldots,\rm{des(\sigma)}$
and labeling the remaining spaces from left to right with
$des(\sigma)+1,\ldots,n$. As an example, let $\sigma=13287546$,
the inv-labeling of $\sigma$ is
\[_8 1_7 3_62_5 8_4 7_3 5_2 4_1 6_0,\]
while the maj-labeling of $\sigma$ is given by
\[_5 1_6 3_4 2_7 8_3 7_2 5_1 4_8 6_0.\]

For $n \geq 2$, we define the map
\[\phi_{inv,n}\colon \{0,1,\ldots,n-1\}\times S_{n-1} \to S_n\]
by setting $\phi_{inv,n}(i,\sigma)$ to be the permutation obtained by inserting
$n$ in the space labeled $i$ in the inv-labeling of $\sigma$.
By changing the inv-labeling to maj-labeling, we
obtain the map $\phi_{maj,n}$. As an example,
 \[\phi_{inv,9}(3,13287546)=132879546 ~~\text{and}~~ \phi_{maj,9}(3,13287546)=132897546. \]

For maps $\phi_{inv,n}$ and $\phi_{maj,n}$, we have the following two lemmas.
\begin{lem}\label{inv}
For $\sigma \in S_{n-1}$ and $i \in \{0,1, \ldots, n-1 \}$, we have \[ inv(\phi_{inv,n}(i,\sigma))=inv(\sigma)+i.\]
\end{lem}
\begin{lem}\label{maj}
For $\sigma \in S_{n-1}$  and $i \in \{0,1, \ldots, n-1 \}$, we have \[ maj(\phi_{maj,n}(i,\sigma))=maj(\sigma)+i.\]
\end{lem}
Lemma \ref{inv} is easy to verified. For a detailed proof of Lemma \ref{maj}, see \cite{Haglund}.

Given a permutation $\pi \in S_n$, let $\pi^{(i)}$ be the
restriction of $\pi$ to the letters $1,2,\ldots,i$ for $1 \leq i \leq n$ .
For $2 \leq i \leq  n$, let
\[(c_i,\pi^{(i-1)})= \phi_{inv,i}^{-1}(\pi^{(i)}),\]  \[(m_i,\pi^{(i-1)})= \phi_{maj,i}^{-1}(\pi^{(i)}).\]
By setting $c_1=0$ and $m_1=0$, we
obtain two sequences $c_1c_2 \cdots c_n$ and $m_1m_2 \cdots m_n$ in $E_n$, where
\[E_n=\{w=w_1 \cdots w_n| w_i\in [0,i-1], 1 \leq i \leq n\}.\]
Define $\gamma(\pi)=c_1c_2 \cdots c_n $ and
$\mu(\pi)=m_1m_2 \cdots m_n $. It is not hard to check that
both $\gamma$ and $\mu$ are bijections. The sequence $c_1c_2 \cdots c_n$ is called the inversion table of $\pi$, while $m_1m_2 \cdots m_n$ is called the
major index table of $\pi$. Moreover, we have
$\sum_{i=1}^{n}c_i=inv(\pi)$ and
$\sum_{i=1}^{n}m_i=maj(\pi).$
As an example, let $\pi=13287546$, then $\pi^{(8)}=\pi$ and
$\mu(\pi)=00204056$ as computed in Table \ref{majtable}.
\begin{table}[h]
\begin{center}
\begin{tabular}{|c|c c|}
  \hline
$i$ &$\pi^{(i-1)}$ & $m_i$ \\\hline
  $8$ & ~~~~~~ $_4 1_5 3_3 2_6  7_2 5_1 4_7 6_0$~~~~~ & $6$\\
  $7$ & $_3 1_4 3_2 2_5 5_1 4_6 6_0$ & $5$ \\
  $6$ & $_3 1_4 3_2 2_5 5_1 4_0$ & $0$ \\
  $5$ & $_2 1_3 3_1 2_4 4_0 $ & $4$ \\
  $4$ & $_2 1_3 3_1 2_0 $& $0$ \\
  $3$ & $_1 1_2 2_0$ & $2$ \\
  $2$ & $_1  1_0$ & $0$\\
  \hline
\end{tabular}
\caption{The computation of the major index table of $13287546$.}
\label{majtable}
\end{center}
\end{table}

Now, we can define the bijection $\varphi$ that maps the inversion number to the major index  by letting $\varphi=\mu^{-1}\gamma$. Clearly, $\varphi$ is a bijection
which proves that statistics $inv$ and $maj$ are equidistributed over $S_n$.

In the following of this section, we will construct a bijection $\rho$ to prove Conjecture \ref{stat}, which is, to some extent, an analogue of the above bijection $\varphi$.

First, we define a stat-labeling of $\sigma \in S_n$.
Label the descents of $\sigma$ and the space after $\sigma_n$  by $0,1,\ldots,des(\sigma)$ from left to right. Label the  space before $\sigma_0$ by $des(\sigma)+1$.
The ascents of
$\sigma$ are labeled from right to left by $des(\sigma)+2, \ldots, n$.
As an example, for $\sigma=13287546$, we have the stat-labeling of $\sigma$ as follows
\[_5 1_8 3_0 2_7 8_1 7_2 5_3 4_6 6_4.\]
Based on the stat-labeling,
we define the map  $\phi_{stat,n}$ for $n \geq 2$.
For a permutation $\sigma=\sigma_1 \sigma_2 \ldots \sigma_{n-1}$ and  $0 \leq i \leq n-1$, define
$\phi_{stat,n}(i,\sigma)$
to be the permutation obtained by inserting
$n$ in the position labeled $i$ in the stat-labeling of $\sigma$.
For instance, $\phi_{stat,9}(7,13287546)=132987546$.

It should be noted that unlike the properties of  $\phi_{inv,n}$ and $\phi_{maj,n}$ stated in Lemma \ref{inv}
and Lemma \ref{maj}, we deduce that  $\phi_{stat,n}(i,\sigma)) \neq stat
(\sigma)+i$ for $i=des(\sigma)+1$. For the map $\phi_{stat,n}$, we have the following property.

\begin{lem}\label{statlem}
For $\sigma \in S_{n-1}$ and $i \in \{0,1,\ldots, des(\sigma),des(\sigma)+2, \ldots, n-1 \}$, we have \[stat(\phi_{stat,n}(i,\sigma))=stat
(\sigma)+i.\]
\end{lem}
\pf
First, we recall that $stat=(ac$-$b)+(ba$-$c)+(cb$-$a)+(ba)$. To prove this lemma, we have to consider the changes of the statistic $stat$ brought by inserting $n$ into $\sigma$. Assume that $\sigma=\sigma_1 \sigma_2\cdots \sigma_{n-1}$, there are three cases for us to consider.
\begin{itemize}
  \item Case 1 $\colon$ $n$ is inserted into the space after $\sigma_{n-1}$.

  Write $\pi=\sigma n$.
  Clearly, the insertion of $n$ does not bring new  $ac$-$b$, $cb$-$a$ and
  $ba$ patterns. While $n$ can form new $ba$-$c$ patterns of $\pi$ with the $ba$ patterns of $\sigma$.
   It follows that $stat(\pi)-stat(\sigma)=(ba)\sigma=des(\sigma)$.
  Notice that  the label of the space after $\sigma_{n-1}$ is $des(\sigma)$. Hence, in this case the lemma holds.

  \item Case 2 $\colon$ $n$ is inserted into a descent.

  Suppose that $\tau$ is the permutation obtained by inserting $n$ to the position
  $i$ and $\sigma_i> \sigma_{i+1}$. Moreover, let this descent be the $k$-th descent from left to right.
  We claim that  $stat(\tau) -stat(\sigma)=k-1$.

   The insertion of $n$ forms some new $ac$-$b$ patterns, and the number of
   these new patterns is $|\{j, ~j>i ~\text{and}~ \sigma_j>\sigma_i\}|$.
  Moreover, $n$ forms $k-1$ new $ba$-$c$ patterns with
   the former $k-1$ descents, while it destroys the  $ba$-$c$ patterns of $\sigma$  by the number $|\{j, ~j>i ~\text{and}~ \sigma_j>\sigma_i\}|$.
   It is easy to verify the functions $(cb$-$a)$ and $(ba)$
   do not change.
   Hence, we conclude that
   $stat(\tau) -stat(\sigma)=k-1$. The claim is verified.
   Notice that the label of the $k$-th descent from left to right is also $k-1$.
   It follows that  the lemma holds for this case.

  \item  Case 3 $\colon$ $n$ is inserted into an ascent.

  Suppose that $p$ is the permutation obtained by inserting $n$ into the position $i$, where $\sigma_i < \sigma_{i+1}$.
  Moreover, we assume that there are $k$ descents to the left of position $i$.
  We claim that $stat(p)-stat(\sigma)=k+n-i+1$.

  Now we proceed to prove this claim.
  The insertion of $n$ to $\sigma$ brings new $ac$-$b$ patterns, the number of which is  $|\{j, ~j \geq i+1 ~\text{and}~ \sigma_j \geq \sigma_{i+1}\}|$.
  Moreover, the insertion of $n$ brings $k$ new $ba$-$c$ patterns with
  the former $k$ descents.
  Also, $n\sigma_{i+1}\sigma_l$, where $l>i+1$ and $\sigma_l<\sigma_{i+1}$, forms a new $cb$-$a$ pattern of $p$.
  Notice that $(ba)p-(ba)\sigma=1$.
  By combining all above,
  we see that  $stat(\tau) -stat(\sigma)=k+n-i+1$. The claim is verified. By the stat-labeling of $\sigma$, we see that the label of this position
  is also $k+n-i+1$. Hence, in this case the lemma holds.
\end{itemize}
By combining the three cases above, we compete the proof.\qed

Based on the stat-labeling,
we define the stat table of a permutation.
Given $\pi \in S_n$,
 for $2 \leq i \leq  n$, let
\[(s_i,\pi^{(i-1)})= \phi_{stat,i}^{-1}(\pi^{(i)}).\]
Set $s_1=0$ and $\nu(\pi)=s_1 s_2 \cdots s_n$. It is easily checked that $\nu$ is a bijection from $S_n$ to $E_n$.
As an example, $\nu(52718346)=01112216$, which is computed in Table \ref{stattable}.
\begin{table}[h]
\begin{center}
\begin{tabular}{|c|c c|}
  \hline
$i$ &$\pi^{(i-1)}$ & $s_i$ \\\hline
  $8$ & ~~~~~~ $_3 5_0 2_7 7_1 1_6 3_5 4_4 6_2 $~~~~~ & $6$\\
  $7$ & $_3 5_0 2_1 1_6 3_5 4_4 6_2 $ & $1$ \\
  $6$ & $_3 5_0 2_1 1_5 3_4 4_2 $ & $2$ \\
  $5$ & $_2 2_0 1_4 3_3 4_1$ & $2$ \\
  $4$ & $_2 2_0 1_3 3_1 $& $1$ \\
  $3$ & $_2 2_0 1_1 $ & $1$ \\
  $2$ & $_1  1_0$ & $1$\\
  \hline
\end{tabular}
\caption{The computation of the stat table of $52718346$.}
\label{stattable}
\end{center}
\end{table}

Now we can give
the definition of the map $\rho$ which proves Conjecture \ref{stat}.
Given $\pi \in S_n$, let $\sigma= \rho(\pi)$, where $\sigma$ can be constructed as follows.
Assume that $F(\pi)=k$ and
$\pi^{(k)}=k p_2 \cdots p_k $. Then, let $\sigma^{(k)}=k (k-p_k) (k-p_{k-1}) \cdots (k-p_2)$. Assume that $s=s_1 s_2 \cdots s_n=\nu(\pi)$. For $k+1 \leq i \leq n$, let $\sigma^{(i)}=\phi_{maj,i}(s_i,\sigma^{(i-1)})$.
Clearly, $\sigma=\sigma^{(n)}$ can be constructed by the procedure above. As an example, we
compute $\rho(\pi)$, where $\pi=52718346$. It is straightforward to see that $k=5$ and $\pi^{(5)}=52134$.
\begin{table}[h]
\begin{center}
\begin{tabular}{|c|c c|}
  \hline
  $i$ & $s_i$ & $\sigma^{(i)}$ \\\hline
  $5$ & $2$ & $_3 5_2 1_4 2_5 4_1 3_0 $ \\
  $6$ & $2$ & $_3 5_4 6_2 1_5 2_6 4_1 3_0$ \\
  $7$ & $1$ & $_3 5_4 6_2 1_5 2_6 4_7 7_1 3_0$\\
  $8$ & $6$ & $56128473$\\
  \hline
\end{tabular}
\caption{The computation of $\rho(52718346)$.}
\label{rho}
\end{center}
\end{table}
By the computation in Table \ref{rho}, we see that
$\rho(52718346)=56128473$. Notice that in this example,
the descent number and the first letter of both of the preimage and image of $\rho$
are the same. In fact, these properties always holds.

\begin{lem}\label{F}
For $\pi \in S_n$, we have $F(\pi)=F(\rho(\pi))$ and
$des(\pi)=des(\rho(\pi))$.
\end{lem}

\pf
Suppose that $\sigma=\rho(\pi)$ and $k=F(\pi)$. We proceed to show that $F(\sigma^{(i)})=F(\pi^{(i)})=k$ and $des(\sigma^{(i)})=des(\pi^{(i)})$ for $k \leq i \leq n$ by induction. Let $\pi^{(k)}=k p_2 \cdots p_k $, then
we have $\sigma^{(k)}=k (k-p_k) (k-p_{k-1}) \cdots (k-p_2)$.
It is routine to check that $F(\sigma^{(k)})=F(\pi^{(k)})=k$ and $des(\sigma^{(k)})=des(\pi^{(k)})$, hence, we omit the details here.

Now assume that $des(\sigma^{(l)})=des(\pi^{(l)})=d$ and
 $F(\sigma^{(l)})=F(\pi^{(l)})=k$ for $k \leq l \leq n-1 $. We proceed to show that
$des(\sigma^{(l+1)})=des(\pi^{(l+1)})$ and $F(\sigma^{(l+1)})=F(\pi^{(l+1)})=k$.

Let $s=s_1 s_2 \cdots s_n=\nu(\pi)$.
By the constructions of $\nu$ and $\rho$, we may see that
$\pi^{(l+1)}=\phi_{stat,l+1}(s_{l+1},\pi^{(l)})$ and
$\sigma^{(l+1)}=\phi_{maj,l+1}(s_{l+1},\sigma^{(l)})$.
Notice that both in the maj-labeling of $\sigma^{(l)}$ and the stat-labeling of $\pi^{(l)}$, the descents and the space
after the last letter are labeled by $\{0,1,\ldots,d\}$, the
space before the first element is labeled by $d+1$, and
the ascents are labeled by $\{d+2, \ldots, l\}$. Since $F(\pi^{(l)})=k=F(\pi)$, it is easy to see that $s_{l+1} \neq d+1$.
It follows that $F(\sigma^{(l+1)})=F(\pi^{(l+1)})=k$.

If $s_{l+1}<d+1$, then $l+1$ is inserted into the descents or the space after the last element
of $\pi^{(l)}$ and $\sigma^{(l)}$. Hence, we deduce that
 $des(\sigma^{(l+1)})=des(\pi^{(l+1)})=d$.
   If $s_{l+1}> d+1$, then $l+1$ is inserted into the ascents
of $\pi^{(l)}$ and $\sigma^{(l)}$. Hence, we deduce that
 $des(\sigma^{(l+1)})=des(\pi^{(l+1)})=d+1$. Combining the
 two cases above, we have $des(\sigma^{(l+1)})=des(\pi^{(l+1)})$.
 Notice that $\sigma^{(n)}=\sigma$ and $\pi^{(n)}=\pi$, we complete  the proof. \qed


Base on the construction of $\rho$ and  Lemma \ref{F}, we have the following theorem.
\begin{thm}\label{bijection}
 The map $\rho$ is an involution on $S_n$.
\end{thm}
\pf
Given a permutation $\pi \in S_n$, it suffices for us
to show that $\rho^2(\pi)=\pi$. That is, writing $\sigma=\rho(\pi)$, we need to show that $\rho(\sigma)=\pi$.

Let $F(\pi)=k$, then
by Lemma \ref{F},
we know that $k=F(\pi)=F(\sigma)$. Write $\pi^{(k)}=k p_2 \cdots p_k$, then we have $\sigma^{(k)}=k (k-p_k) \cdots (k-p_2) $. Assume that $\nu(\pi)=s=s_1 s_2 \cdots s_n$ and
$\mu(\sigma)=m_1 m_2 \cdots m_n$.
By the construction of $\rho$, we have
$\sigma^{(i)}=\phi_{maj,i}(s_i,\sigma^{(i-1)})$ for $k+1 \leq i \leq n$. It follows that $m_i=s_i$ for $k+1 \leq i \leq n$.

Suppose that $\alpha=\rho(\sigma)$, in the following,
we proceed to show that $\alpha^{(i)}=\pi^{(i)}$ for $k \leq i \leq n $ by induction. By the definition of $\rho$, we have
$\alpha^{(k)}=kp_2\cdots p_k=\pi^{(k)}$. Assume that
$\alpha^{(i)}=\pi^{(i)}$ holds for $k \leq i \leq n-1$,
we aim to show that $\alpha^{(i+1)}=\pi^{(i+1)}$. To achieve
this, we need to mention the following property of the maj-labeling and the stat-labeling of a single permutation.

For a permutation $p\in S_n$, assume the maj-labeling of $p$ is $f_0f_1 \cdots f_n$,
while the stat-labeling of $p$ is $h_0h_1 \cdots h_n$.
Then it is easily checked that

\begin{equation}\label{label}
f_i+h_i = \left\{
  \begin{array}{ll}
des(p), & \mbox{if $i$ is a descent of $p$ or $i=n$,}\\[3pt]
n+des(p)+2, & \mbox{if $i$ is an ascent of $p$,}\\[6pt]
2des(p)+2, & \mbox{if $i=0$.}
 \end{array} \right.
\end{equation}

Let $\nu(\sigma)=l_1 l_2 \cdots l_n$ and $d=des(\alpha^{(i)})$.
Then, we have
\begin{equation*}
\alpha^{(i+1)} =\phi_{maj,i+1}(l_{i+1},\alpha^{(i)}) = \left\{
  \begin{array}{ll}
\phi_{stat,i+1}(d-l_{i+1},\alpha^{(i)}), & \mbox{if $0 \leq l_{i+1} \leq d $,}\\[3pt]
\phi_{stat,i+1}(n+d+2-l_{i+1},\alpha^{(i)}), & \mbox{if $d+2 \leq l_{i+1} \leq n $.}
 \end{array} \right.
\end{equation*}
By the proof of Lemma \ref{F}, we see that $des(\sigma^{(i)})=des(\alpha^{(i)})=d$.
Recall that $\nu(\sigma)=l_1 l_2 \cdots l_n$ and $\mu(\sigma)=m_1 m_2 \cdots m_n$. Hence, it
follows from (\ref{label}) that
\begin{align*}
\alpha^{(i+1)}=& ~\phi_{stat,i+1}(m_{i+1},\alpha^{(i)})\\[3pt]
    =& ~\phi_{stat,i+1}(s_{i+1},\pi^{(i)})\\[3pt]
    =& ~\pi^{(i+1)}.
\end{align*}
Notice that $\alpha=\alpha^{(n)}$ and $\pi=\pi^{(n)}$.
Hence, we have $\pi=\rho(\sigma)$, namely, $\rho^2(\pi)=\pi$.
This completes the proof.\qed

Indeed, the involution
$\rho$ also preserves some addtional statistics, which is
stated in the following proposition.
\begin{prop}\label{stat=maj}
For any $\pi=\pi_1 \pi_2 \cdots \pi_n \in S_n$, we have  $maj(\rho(\pi))=stat(\pi)$ and $stat(\rho(\pi))=maj(\pi)$.
\end{prop}
\pf
Write $\sigma=\rho(\pi)$. Let $k=F(\pi)$ and
$\pi^{(k)}=k p_2 \cdots p_k$. Then we have $\sigma^{(k)}=
k (k-p_k) \cdots (k-p_2)$. In the following, we will show
that $maj(\sigma^{(i)})=stat(\pi^{(i)})$ for $k \leq i \leq n$ by induction.

First, we show that $maj(\sigma^{(k)})=stat(\pi^{(k)})$.
By the proof of Lemma \ref{F}, we have $des(\pi^{(k)})=des(\sigma^{(k)}).$
Hence $(ba)\pi^{(k)}=(ba)\sigma^{(k)}$.  It suffices for us to show that
  \begin{equation}\label{**}
    (ac-b)\pi^{(k)}+(ba-c)\pi^{(k)}+(cb-a)\pi^{(k)}=
  (a-cb)\sigma^{(k)}+(b-ca)\sigma^{(k)}+(c-ba)\sigma^{(k)}
  \end{equation}
  Given a pattern $p=p_1 p_2 \cdots p_k$, by putting a line under $p_1$(resp. $p_k$), we mean that an instance of $p$ must begin(resp. end) with the leftmost(resp. rightmost) letter
  of the permutation. By putting a dot under $p_1$, we mean that
  an instance of $p$ must not begin with the leftmost letter.
  For instance, $(\b{c}b$-$a)$ is a function which maps a permutation, say $\tau=\tau_1 \tau_2 \cdots \tau_n$, to
  $|\{i, \tau_1 \tau_2\tau_{i} ~\text{forms a 321 pattern of}~\tau\}|.$

  By the definition of $\rho$, we know that
  \begin{align*}
    &(ac-b)\pi^{(k)}+(ba-c)\pi^{(k)}+(cb-a)\pi^{(k)}\\
    =&(ac-b)\pi^{(k)}+(ba-c)\pi^{(k)}+(\b{c}b-a)\pi^{(k)}+(\d{c}b-a)\pi^{(k)}\\
    =&(b-ac)\sigma^{(k)} + (a-cb)\sigma^{(k)} + (\b{c}-b-\b{a})\sigma^{(k)} + (\d{c}-ba)\sigma^{(k)}
  \end{align*}
   Hence, to prove (\ref{**}), we need to prove that for any permutation $p \in S_k$ with $p_1=k$,
   \begin{equation} \label{***}
    (\b{c}-ba)p+(b-ca)p=(\b{c}-b-\b{a})p+(b-ac)p.
   \end{equation}
  We define sets $A(p),C(p)$ and multisets $ B(p),D(p)$ as follows.
   \begin{align*}
    A(p)&=\{p_i \colon p_1p_i p_{i+1} \text{ forms a $321$ pattern of }   p\},\\
    B(p)&=\{p_i \colon p_ip_j p_{j+1} \text{ forms a $231$ pattern of }   p\},\\
    C(p)&=\{p_i \colon p_1p_i p_{m} \text{ forms a $321$ pattern of }   p\},\\
    D(p)&=\{p_i \colon p_ip_j p_{j+1} \text{ forms a $213$ pattern of }   p\}.
   \end{align*}
 To prove (\ref{***}), it is enough to show that $A \cup B =C\cup D$,
 where the union operator is a multiset union.
 First, we show that
$C\cup D  \subseteq A \cup B$.

Let $p_k=a$, then we know that $C(p)=\{a+1, a+2, \ldots, k-1\}$.
If $p_i \in C(p)$ is a descent top, it is easy to see that
$p_i \in A(p)$.
 If $p_i \in C(p)$ is an ascent bottom,
we claim that $p_i \in B(p)$. This claim will be
proved together with case 4 in the following.

Given an element $p_i$ in $D(p)$, if the multiplicity of $p_i$
is $x$,  there exists a set \[\{j_1,j_1+1, j_2, j_2+1,\ldots,j_x,j_x+1\},\] which is ordered by increasing
order, satisfying that \[p_ip_{j_1}p_{j_1+1},~ p_ip_{j_2}p_{j_2+1},~\ldots,~
         p_ip_{j_x}p_{j_x+1}\] are instances of $213$ pattern.
         We claim that there exists $j_1+1 \leq r_1 < j_{2}$ such that
         $p_i p_{r_1}p_{r_1+1}$ forms a $231$ pattern.

         Choose the smallest $g_1$ such that $j_1+1 \leq g_1< j_2$ and
         $g_1$ is a descent. If $p_{g_1+1}<p_i$, then $p_i p_{g_1}p_{g_1+1}$ forms a $231$ pattern, the claim is verified.
         Otherwise, we seek  the smallest $g_1+1 \leq g_2< j_2$ such that
         $g_2$ is a descent. If $p_{g_2+1}<p_i$, the claim is verified.
         If not, we repeat the above process. Since $p_i > p_{j_2}$,
         the process must be terminated. Hence, the claim is verified.
          By a similar means, we deduce that there exists $j_l+1 \leq r_l < j_{l+1}$ where $1 \leq l \leq x-1$ such that $p_i p_{r_l}p_{r_l+1}$ forms a $231$ pattern. The claim is verified.

To analyze the element $p_i$ in $D(p)$, we consider four cases.
\begin{itemize}
  \item  $p_i$ is a descent top and $1 < p_i \leq a-1$.

         Suppose that  the multiplicity of $p_i$ of this type in $D(p)$ is $x$.
          By the above statement, we see that $p_i p_{r_l}p_{r_l+1}$ forms a $231$ pattern, where $j_l+1 \leq r_l < j_{l+1}$ and $1 \leq l \leq x-1$. Thus, we deduce that there are $x-1$ $p_i$s in $B(p)$.
         Notice that there is one $p_i$ left in $B(p)$. Clearly,
         we can set this $p_i$ to be the element of $A(p)$
         consisting of $p_1$, $p_i$ and $p_{i+1}$.

  \item  $p_i$ is an ascent bottom and $1 \leq p_i \leq a-1$.

       Suppose that the multiplicity of $p_i$ of this type in $D(p)$ is $x$.
       Similarly, we know that $p_i p_{r_l}p_{r_l+1}$ forms a $231$ pattern, where $j_l+1 \leq r_l < j_{l+1}$ and $1 \leq l \leq x-1$.
       Since $p_i < p_{i+1}$, we have $j_1>i+1$. Based on this,
       it can be easily seen that there exists $r_0$ such that $i+1 \leq r_0 < j_1$ and
       $p_i p_{r_0}p_{r_0+1}$ forms a $231$ pattern.
       Hence, we deduce that in this case there are $x$ $p_i$s in $B(p)$.

       \item  $p_i$ is a descent top and $a+1 \leq p_i \leq k-1$.

       Suppose that  the multiplicity of $p_i$ of this type in $D(p)$ is $x$.
       Similarly, we deduce that $p_i p_{r_l}p_{r_l+1}$ forms a $231$ pattern, where $j_l+1 \leq r_l < j_{l+1}$ and $1 \leq l \leq x-1$.
       What's more, it follows from  $p_i>a$  that there exists $j_x \leq r_x<m$ such that
       $p_i p_{r_{x}}a$ forms a $231$ pattern.
       Hence, we deduce that $p_i$ in $B(p)$ and its multiplicity
       is $x$.

       \item  $p_i$ is an ascent bottom and $a+1 \leq p_i \leq k-1$.

       Suppose that  the multiplicity of $p_i$ of this case in $D(p)$ is $x$. Notice that $p_i$ is also an element of $C(p)$ with multiplicity equals $1$. Hence, in this case, we have to prove that
       there are $x+1$ $p_i$s in $B(p)$.

       Similarly with the above cases, we deduce that $p_i p_{r_l}p_{r_l+1}$ forms a $231$ pattern, where $j_l+1 \leq r_l < j_{l+1}$ and $1 \leq l \leq x-1$.
       Since $p_i>a$, there exists $j_x \leq r_x<m$ such that
       $p_i p_{r_{x}}a$ forms a $231$ pattern.
       By $p_i < p_{i+1}$, we have $j_1>i+1$. Based on this,
       it can be easy seen that there exists $r_0$ such that $i+1 \leq r_0 < j_1$ and
       $p_i p_{r_0}p_{r_0+1}$ forms a $231$ pattern.
       Hence, we deduce that $p_i$ in $B(p)$ and its multiplicity
       is $x+1$.
\end{itemize}

        Combining all above, we deduce that $C\cup D  \subseteq A \cup B$.
        By a similar analysis, we can prove that $A \cup B  \subseteq C\cup D$.
        We omit it here.  As an example, if $p=978452613$, we have
       $A(p)=\{5,6,8\}$, $B(p)=\{2,4,4,5,7\}$, $C(p)=\{4,5,6,7,8\}$
       and $D(p)=\{2,4,5\}$. It can be verified that
       $A \cup B =C\cup D$. This proves that $maj(\sigma^{(k)})=stat(\pi^{(k)})$.

Now assume that $maj(\sigma^{(i)})=stat(\pi^{(i)})$ for
$k \leq i \leq n-1$, we proceed to show that $maj(\sigma^{(i+1)})=stat(\pi^{(i+1)})$.
Write $\nu(\pi)=s_1 s_2 \cdots s_n$.
Then, by Lemma \ref{maj}, Lemma \ref{statlem} and the
construction of $\rho$, we have
\begin{align*}
maj(\sigma^{(i+1)})=&~ s_{i+1}+maj(\sigma^{(i)}) \\[3pt]
    =& ~s_{i+1}+stat(\pi^{(i)})  \\[3pt]
    =& ~ stat(\pi^{(i+1)}).
\end{align*}
This proves $maj(\sigma^{(i)})=stat(\pi^{(i)})$ for $k \leq i \leq n$.
Notice that $\pi=\pi^{(n)}$ and $\sigma=\sigma^{(n)}$.
We deduce that $maj(\sigma)=maj(\rho(\pi))=stat(\pi)$.
By Theorem \ref{bijection}, we see that $\rho$ is an involution. This implies that $stat(\rho(\pi))=maj(\pi)$.
This completes the proof.
\qed

Combining Lemma \ref{F}, Theorem \ref{bijection} and Proposition \ref{stat=maj}, we give a
proof of Conjecture \ref{stat}.

It should be mentioned that Burstein \cite{Burstein} provided a direct bijective proof of a refinement of Conjecture \ref{stat}.
The bijection $\chi$ is given as follows.
Given a permutation $\pi \in S_n$ with $F(\pi)=k$. Let $\pi'=\chi(\pi)$ with $\pi'(1)=k$ and
\begin{equation}\label{label}
\pi'(i) = \left\{
  \begin{array}{ll}
k-\pi(n+2-i), & \mbox{if $\pi(n+2-i)<k$,}\\[3pt]
n+k+1-\pi(n+2-i), & \mbox{if if $\pi(n+2-i)>k$.}
 \end{array} \right.
\end{equation}
In addition to preserving the statistics $des$ and $F$, the bijection $\chi$
also preserves the statistic $adj$, while our bijection does not.
As an example, set $\pi=543617982$, then $\rho(\pi)=\sigma=539784621$ and
 $\chi(\pi)=\pi'=537684921$. It can be checked that
 $adj(\pi) = adj(\pi')$, while
 $adj(\pi) \neq adj(\sigma)$. Moreover,
 it is easily seen that for $\pi \in S_n$ with $F(\pi)=n$,
 we have $\rho(\pi)=\chi(\pi)=\sigma$. We note that
 in this case both Burstein and us have to prove $maj(\sigma)=stat(\pi)$. Different form our proof
 in Proposition \ref{stat=maj}, Burstein gave the following
 two relations, which implies that $maj(\sigma)=stat(\pi)$.
 \[maj(\pi)+stat(\pi)=(n+1)des(\pi)-(F(\pi)-1),\]
 \[maj(\pi)+maj(\sigma)=(n+1)des(\pi)-(F(\pi)-1).\]

\vspace{0.5cm}
\noindent{\bf Acknowledgments.}  We wish to thank the
 anonymous referees for their valuable comments and suggestions.

\end{document}